\documentclass[12pt]{article}
\usepackage{latexsym}
\usepackage{amsmath,amssymb}
\usepackage{amsthm}
\usepackage{amsfonts}
\begin{document}
\abovedisplayskip=10pt plus 1pt minus 3pt \belowdisplayskip=10pt
plus 1pt minus 3pt
\def\bq{\begin{equation}\hspace*{1.0cm}}
\def\eq{\end{equation}}
\def\ba{\begin{array}{lllll}}
\def\br{\begin{array}{rlllll}}
\def\ea{\end{array}}
\def\ST{\songti\rm\relax}   \def\HT{\heiti\bf\relax}
\def\FS{\fangsong\rm\relax} \def\KS{\kaishu\rm\relax}
\def\dx{{\rm d}x}\def\dt{{\rm d}t}
\def\dpr{\displaystyle\prod} \def\dmin{\displaystyle\min}
\def\dlim{\displaystyle\lim} \def\dmax{\displaystyle\max}
\def\dcup{\displaystyle\bigcup} \def\dcap{\displaystyle\bigcap}
\def\dsup{\displaystyle\sup} \def\dinf{\displaystyle\inf}
\def\dint{\displaystyle\int} \def\dsum{\displaystyle\sum}
\def\va{\varepsilon} \def\la{\lambda} \def\ir{I\kern-.25em R}
\def\bv{{\bf {v}}} \def\bw{{\bf{w}}} \def\bk{{\bf{k}}}\def\bn{{\bf{n}}}
\def\bx{{\bf{x}}}
\def\R{\mathbb{R}}
\def\oiint{{\displaystyle\bigcirc\kern-3.2ex\int\kern-2ex\int}}
\def\pmb#1{\mbox{\boldmath $#1$}} \def\text#1{\hbox{#1}}
\def\MCH#1#2{\setbox0=\hbox{\raise#1\hbox{#2}}\smash{\box0}}
\def\SEC#1{\vspace{0.5\baselineskip}\vbox{\vspace*{1\baselineskip}
\centering\large\HT\zihao{-5}#1}\vspace*{1\baselineskip}\rm}
\def\REF#1{\par\hangindent\parindent\indent\llap{#1\enspace}\ignorespaces}
\def\Cal#1{{\cal #1}}
\renewcommand{\theequation}{\thesection .\arabic{equation}}
\def\ST{\songti\rm\relax}   \def\HT{\heiti\bf\relax}
\def\FS{\fangsong\rm\relax} \def\KS{\kaishu\rm\relax}
\def\dpr{\displaystyle\prod} \def\dmin{\displaystyle\min}
\def\dlim{\displaystyle\lim} \def\dmax{\displaystyle\max}
\def\dcup{\displaystyle\bigcup} \def\dcap{\displaystyle\bigcap}
\def\dsup{\displaystyle\sup} \def\dinf{\displaystyle\inf}
\def\dint{\displaystyle\int} \def\dsum{\displaystyle\sum}
\def\va{\varepsilon} \def\la{\lambda} \def\ir{I\kern-.25em R}
\def\oiint{{\displaystyle\bigcirc\kern-3.2ex\int\kern-2ex\int}}
\def\pmb#1{\mbox{\boldmath $#1$}} \def\text#1{\hbox{#1}}
\def\MCH#1#2{\setbox0=\hbox{\raise#1\hbox{#2}}\smash{\box0}}
\def\SEC#1{\vspace{0.5\baselineskip}\vbox{\vspace*{3.0\baselineskip}
\centering\large\HT\zihao{-4}#1}\vspace*{3.0\baselineskip}\rm}
\def\REF#1{\par\hangindent\parindent\indent\llap{#1\enspace}\ignorespaces}
\def\Cal#1{{\cal #1}}
\renewcommand{\theequation}{\thesection.\arabic{equation}}
\def\disp{\displaystyle}

\newpage
\thispagestyle{empty}
\begin{center}
\renewcommand{\thefootnote}{\fnsymbol{footnote}}
{\large \bf Local well-posedness and  blow up criterion for
the\\[3mm]
Inviscid Boussinesq system in  H\"{o}lder spaces}
\end{center}

\begin{center}
{\normalsize{Xiaona CUI }\footnote{ e-mail:
 cxn009@126.com}}

{\normalsize\it School of Mathematical Sciences of Graduate
University of Chinese Academy of Sciences, Beijing 100049, PRC}

{\normalsize{Changsheng DOU }\footnote{ e-mail:
 douchangsheng@163.com}}

{\normalsize\it School of Mathematical Sciences,

Capital Normal University, Beijing 100048,PRC}

{\normalsize{Quansen JIU }\footnote{The research is partially
supported by National Natural Sciences Foundation of China (No.
10431060).  e-mail:
 jiuqs@mail.cnu.edu.cn}}

{\normalsize\it School of Mathematical Sciences,

Capital Normal University, Beijing 100048,PRC}

\end{center}

\begin{center}
\parbox{0.8\hsize}{\parindent=4mm
{\it Abstract}:\ We prove the local in time existence and a blow up
criterion of solution in the H\"{o}lder spaces for the inviscid
Boussinesq system in $R^{N},N\geq2$, under the assumptions that the
initial values $\theta_{0},u_{0}\in C^{r}$, with $r>1$.

 {\it Key
Words}: inviscid Boussinesq system, local well-posedness, blow-up
criterion

{\it AMS(2000)Subject Classification}: 35Q35, 76B03}

\end{center}

\vspace{0.5cm}

\setcounter{section}{1} \thispagestyle{empty} \pagestyle{myheadings}
\setcounter{page}{1}

 \setcounter{equation}{0}

{\large \bf 1. Introduction }

\vspace{3mm}

   The Cauchy problem for the Boussinesq system in
$R^{N} (N\ge 2)$ can be written as
\begin{equation}\label{1.1}
 \left\{
 \begin{array}{llll}
\partial_{t}\theta+u\cdot\nabla\theta-\kappa\Delta\theta=g,
\\
\partial_{t}u+u\cdot\nabla u+\nabla\Pi-\nu\Delta u=f,
 \\
 {\rm div} u=0,
\end{array}\right.
 \end{equation}
with  initial data
\begin{equation}
 \theta|_{t=0}=\theta_0,\ \ \  u|_{t=0}=u_0.
\end{equation}
Here  $u(x,t), (x,t)\in R^{N}\times(0,\infty), N\ge 2,$ is the
velocity vector field, $\theta(x,t)$ is the scalar temperature,
$\Pi(x,t)$ is the scalar pressure, $f(x,t)$ is the external forces,
which is a vector function, and $g$ is a known scalar function.
$\nu\geq0$ is the kinematic viscosity, and $\kappa\geq0$ is the
thermal diffusivity.

The Boussinesq system is extensively used   in the atmospheric
sciences and oceanographic turbulence in which rotation and
stratification are important (see [18, 19] and references therein).
When $\kappa>0$ and $\nu>0$, the 2D Boussinesq system with $g=0,
f=\theta e_{2} (e_2=(0,1))$ has been well-understood (see [3], [14],
[1]). When $\kappa>0$ and $\nu=0$ or $\kappa=0$ and $\nu>0$, the
Boussinesq system is usually called the partial viscosity one. In
these cases, the Boussinesq system has also been extensively and
successfully studied. In particular, in the case $\kappa>0$ and
$\nu=0$, D. Chae ([4]) proved that the Boussinesq system is globally
well-posedness in $R^m$ for any $m\ge 3$ and this result was
extended by  T. Hmidi and S. Keraani ([15]), R. Danchin and M. Paicu
([11])  to rough initial data in Besov space framework. In the case
$\kappa=0$ and $\nu>0$, the global well-posedness was proved by D.
Chae ([4]), T. Y. Hou and C. Li ([16]) in $H^m(R^2))$ space with
$m\ge 3$.

When $\kappa=0$ and $\nu=0$,  the  Boussinesq system (1.1) becomes
the inviscid one.  In this case, it is clear that if $\theta\equiv
0$, the inviscid Boussinesq system reduces to the classical Euler
equations.  And the two-dimensional Boussinesq system can be used a
model for the three-dimensional axisymmetric Euler equations with
swirl (see [12]). However,  the global well-posedness problem of the
inviscid Boussinesq system is still completely open in general (an
exceptional case is the two-dimensional Euler equations which
correspond to $\theta\equiv 0$, see [8] and references therein).
Local existence and blow-up criteria have been  established for the
inviscid Boussinesq system (see [6, 7, 13, 20] and references
therein). In particular, D. Chae, S. K. Kim and H. S. Nam considered
in [6] the inviscid Boussinesq system  with $g=0$ and $f=\theta
f_{1}$, where $f_{1}$ satisfies that $curl f=0$ and $f\in
L_{loc}^{\infty}([0,\infty);W^{1,\infty}(R^{2}))$. They proved that
there exists a unique  and local $C^{1+\gamma}$ solution of the
inviscid  Boussinesq system with initial data $u_{0},\theta_{0}\in
C^{1+\gamma},\omega_{0},\Delta\theta_{0}\in L^{q}$ for $0<\gamma<1$
and $1<q<2$, where $\omega_0$ is the initial vorticity of the
initial velocity $u_0$.  They also proved a blow-up criterion for
the local solution, which says that the gradient of the passive
scalar $\theta$ controls the breakdown of $C^{1+\gamma}$ solutions
of the Boussinesq system.

In this paper, we devote to the following inviscid  Boussinesq
system:
\begin{equation}
 \left\{
 \begin{array}{llll}
\partial_{t}\theta+u\cdot\nabla\theta=0,
\\
\partial_{t}u+u\cdot\nabla u+\nabla\Pi=\theta e_{N},
 \\
 {\rm div} u=0,
\end{array}\right.
 \end{equation}
with the initial data
\begin{equation}
 \theta|_{t=0}=\theta_0,\ \ \  u|_{t=0}=u_0.
\end{equation}

In [8], the local in time existence and  blow up criterion for the
Euler equations with the initial data $u_{0}\in C^{r} (r>1)$ was
proved. Here we will extend the approaches and results  in [8] to
the inviscid Boussinesq system (1.3). It should be noted that our
results will relax the initial conditions in [6]. More precisely, we
prove that the inviscid Boussinesq system has a local and unique
$C^{r} (r>1)$-solution under the assumptions that the initial data
$u_{0},\theta_{0}\in C^{r} (r>1)$. It does not require that
$\omega_{0},\Delta\theta_{0}\in L^{q}$ for some $1<q<2$ which are
needed in [6].

The plan of the paper is as follows. In Section 2, we give some
preliminaries and our main results.  To prove our main results, we
first present a priori estimates in Section 3 and then in Section 4
we construct approximate solutions and furthermore prove that they
are Cauchy sequences in  appropriate  H\"older spaces. Lastly, in
Section 5, we give the proof of the blow-up criterion.

\vspace{0.2cm}

\renewcommand{\thefootnote}{\fnsymbol{footnote}}
\thispagestyle{empty} \setcounter{section}{2}
\setcounter{subsection}{0}
 \pagestyle{myheadings}

\setcounter{equation}{0}
 {\large \bf 2. Preliminaries and
Main Results }

\vspace{3mm}

Let us start with the definition of the dyadic decomposition of the
full space $R^{N}$ (see[8]).

\textbf{Proposition 2.1}\quad Denote by ${\mathcal{C}}$ the annulus
of centre 0, short radius 3/4 and long radius 8/3. Then there exist
two positive radial functions $\chi$and $\phi$ belonging
respectively to $C_{0}^{\infty}(B(0,4/3))$ and
$C^{\infty}_{0}({\mathcal{C}})$ such that
$$
\chi(\xi)+\sum_{q\geq-1}\varphi(2^{-q}\xi)=1,
$$
$$
|p-q|\geq2\Rightarrow \mbox{supp}\varphi(2^{-q}\cdot)\cap
\mbox{supp}\varphi(2^{-p}\cdot)=\emptyset,
$$
$$
q\geq1\Rightarrow \mbox{supp}\chi\cap
\mbox{supp}\varphi(2^{-q}\cdot)=\emptyset.
$$
If $\widetilde{\mathcal {C}}=B(0,2/3)+\mathcal {C}$, then
$\widetilde{\mathcal {C}}$ is an annulus and we have
$$
|p-q|\geq5\Rightarrow 2^{p}\widetilde{\mathcal {C}}\cap2^{q}\mathcal
{C}=\emptyset,
$$
$$
\frac{1}{3}\leq\chi^{2}(\xi)+\sum_{q\geq0}\varphi^{2}(2^{-q}\xi)\leq1.
$$

\textbf{Notation 2.1}\quad For the inhomogeneous Besov spaces, we
have the notations
$$
h={\mathcal {F}}^{-1}\varphi\quad \mbox{and}\quad
\widetilde{h}={\mathcal {F}}^{-1}\varphi,
$$
$$
\Delta_{-1}u=\chi(D)u={\mathcal
{F}}^{-1}(\chi(\xi)\widehat{u}(\xi)),
$$
$$
\mbox{if}\ q\geq0,\ \Delta_{q}u=\varphi(2^{-q}D)u=2^{qN}\int
h(2^{q}y)u(x-y)dy,
$$
$$
\mbox{if}\  q\leq-2,\  \Delta_{q}u=0,
$$
$$
S_{q}u=\sum_{p\leq
q-1}\Delta_{p}u=\chi(2^{-q}D)u=2^{qN}\int\widetilde{h}(2^{q}y)u(x-y)dy.
$$

The product $uv$ can be formally divided into three parts as follows
(see [2]) :
$$
uv=T_{u}v+T_{v}u+R(u,v),
$$
where
$$
T_{u}v=\sum_{q}S_{q-1}\triangle_{q}v,\ \ \ \
R(u,v)=\sum_{q}\triangle_{q}u(\sum^{1}_{j=-1}\triangle_{q+j}v),
$$
$T_{u}v$ is called paraproduct of $v$ by $u$ and $R(u,v)$ the
remainder term.

\textbf {Proposition 2.2}\quad (Bernstein's Inequality)  Let
$(r_{1}, r_{2})$ be a pair of strictly positive numbers such that
$r_{1}<r_{2}$. There exists a constant $C$ such that for every
nonnegative integer $k$, and for every $1\leq a\leq b$ and for all
function $u\in L^{a}(R^{N})$, we have
$$
\mbox{supp}\hat{u}\in B(0,\lambda r_{1})\quad\Rightarrow\quad
\sup_{|\alpha|=k}\|\partial^{\alpha}u\|_{L^{b}}\leq
C^{k}\lambda^{k+N(\frac{1}{a}-\frac{1}{b})}\|u\|_{L^{a}},
$$
$$
\mbox{supp}\hat{u}\in {\mathcal{C}}(0,\lambda r_{1},\lambda
r_{2})\quad\Rightarrow\quad
C^{-k}\lambda^{k}\|u\|_{L^{a}}\leq\sup_{|\alpha|=k}\|\partial^{\alpha}u\|_{L^{a}}\leq
C^{k}\lambda^{k}\|u\|_{L^{a}},
$$
where $\hat{u}$ denotes the Fourier transform of $u$; $B(0,r)$
refers a ball with the center 0 and radius $r$; and
${\mathcal{C}}(0,r_{1},r_{2})$ denotes analogous as before a ring of
center 0, short radius $r_{1}$ and long radius $r_{2}$.

\textbf{Definition 2.1}\quad Let $s\in R,\ p,\ q\in[1,\infty]$. The
inhomogeneous Besov space $B^{s}_{p,q}$ and the homogeneous Besov
space $\dot{B}^{s}_{p,q}$ is defined as a space of $f\in S^{\prime}$
(tempered distributions) such that
$$
\|f\|_{B^{s}_{p,q}}=(\sum_{j=-1}^{\infty}2^{jqs}\|\Delta_{q}f\|^{q}_{L^{p}})^{\frac{1}{q}}<\infty,\quad
\|f\|_{\dot{B}^{s}_{p,q}}=(\sum_{j=-\infty}^{\infty}2^{jqs}\|\dot{\Delta}_{q}f\|^{q}_{L^{p}})^{\frac{1}{q}}<\infty.
$$
In case $q=\infty$, the expressions are understood as
$$\displaystyle{\|f\|_{B^{s}_{p,\infty}}=\sup_{j\geq-1}2^{js}\|\Delta_{j}f\|_{L^{p}},\quad
\|f\|_{\dot{B}^{s}_{p,\infty}}=\sup_{j\in
Z}2^{js}\|\dot{\Delta}_{j}f\|_{L^{p}} }.$$

And, in the case $p=q=\infty$, we write
$B^{s}_{\infty,\infty}\triangleq C^{r}$, and
$\dot{B}^{s}_{\infty,\infty}\triangleq \dot{C}^{r}$. So the norms
are written as
$\|u\|_{\dot{B}^{r}_{\infty,\infty}}\triangleq\|u\|_{\dot{C}^{r}}$
and $\|u\|_{C^{r}}\triangleq\|u\|_{r}$.

\textbf{Lemma 2.1}\quad Let $v$ be a smooth free-divergence vector
field and $f$ a smooth function. Then we have
$$
\|[v\cdot \nabla,\Delta_{q}]f\|_{L^{\infty}}\leq
C(r)2^{-qr}\|f(t)\|_{r}\|\nabla v\|_{L^{\infty}} (q\geq-1,r>1).
$$
The detailed proof of the lemma is referred to [8].

The following embeddings will be used later.

\textbf{Lemma 2.2}

(1) Let $r>\frac{N}{p}+1$, then $ B^{r-1}_{p,q}\hookrightarrow
L^{\infty}\hookrightarrow\dot{B}^{0}_{\infty,\infty}$. When
$p=q=\infty$, we have the special case
$B^{r}_{\infty,\infty}\hookrightarrow L^{\infty}\hookrightarrow
\dot{B}^{0}_{\infty,\infty},\ \mbox{with}\ r>0$.

(2) Let $s>0,\ 1\leq p,q<\infty$,\ and\
$s_{1}-\frac{N}{p_{1}}=s_{2}-\frac{N}{p_{2}},\ s_{1},s_{2}\in R,\
s_{1}<s_{2}$, then we get $B^{s}_{p,q}\hookrightarrow
\dot{B}^{s}_{p,q} \ \mbox{and}\
\dot{B}^{s_{1}}_{p_{1},q}\hookrightarrow \dot{B}^{s_{2}}_{p_{2},q}$.

(3) $B^{r}_{\infty,\infty}\hookrightarrow\dot{B}^{1}_{\infty,1}\
\mbox{with}\ r>1$,\ and\
$B^{r}_{\infty,\infty}\hookrightarrow\dot{B}^{r}_{\infty,\infty}\
\mbox{for}\ r>0$.

\textbf{Proof.} (1) is proved in [5] and (2) is proved in [8] and
[10].

To prove (3), in view of the definition of the norm in the
homogeneous Besov space, we estimate
$$
\begin{array}{rl}
\|u\|_{\dot{B}^{1}_{\infty,1}}&=\displaystyle{\sum_{q\geq0}2^{q}\|\dot{\Delta}_{q}u\|_{L^{\infty}}+\sum_{q\leq-1}2^{q}\|\dot{\Delta}_{q}u\|_{L^{\infty}}}\nonumber\\
&\leq\displaystyle{\|u\|_{B^{1}_{\infty,1}}+(\sum_{q\leq-1}2^{q})\|u\|_{L^{\infty}}}\nonumber\\
&\leq\|u\|_{B^{1}_{\infty,1}}+\|u\|_{L^{\infty}}\nonumber\\
&\lesssim\|u\|_{r}\nonumber,
\end{array}
$$
where we used the embedding $B^{r}_{\infty,\infty}\hookrightarrow
L^{\infty}(r>0)$.

Similarly, we obtain
$$
\begin{array}{rl}
\|u\|_{\dot{B}^{r}_{\infty,\infty}}&=\displaystyle{\sup_{q\in
Z}2^{qr}\|\dot{\Delta}_{q}u\|_{L^{\infty}}}=\displaystyle{\sup_{q\geq0}2^{qr}\|\dot{\Delta}_{q}u\|_{L^{\infty}}+\sup_{q\leq-1}2^{qr}\|\dot{\Delta}_{q}u\|_{L^{\infty}}}\nonumber\\
&\leq\|u\|_{B^{r}_{\infty,\infty}}+\displaystyle{\sup_{q\leq-1}\|\dot{\Delta}_{q}u\|_{L^{\infty}}}\nonumber\\
&=\|u\|_{B^{r}_{\infty,\infty}}+\|u\|_{\dot{B}^{0}_{\infty,\infty}}\nonumber\\
&\lesssim\|u\|_{B^{r}_{\infty,\infty}}\nonumber,
\end{array}
$$
in which we used the embedding $ L^{\infty}\hookrightarrow
\dot{B}^{0}_{\infty,\infty} $.

{\bf Remark 2.1.} As a special case of (2):
$B^{N}_{1,1}\hookrightarrow \dot{B}^{N}_{1,1}$, which we will use in
this paper, the proof is direct:
$$
\begin{array}{rl}
\|u\|_{\dot{B}^{N}_{1,1}}&=\displaystyle{\sum_{q\in
Z}2^{qN}\|\dot{\Delta}_{q}u\|_{L^{1}}}\nonumber\\
&=\displaystyle{\sum_{q\geq0
}2^{qN}\|\dot{\Delta}_{q}u\|_{L^{1}}+\sum_{q\leq-1}2^{qN}\|\dot{\Delta}_{q}u\|_{L^{1}}}\nonumber\\
&\leq\displaystyle{\sum_{q\geq-1}2^{qN}\|\Delta_{q}u\|_{L^{1}}+(\sum_{q\leq-1}2^{qN})2^{N}(2^{-N}\|\Delta_{-1}u\|_{L^{1}})}\nonumber\\
&\lesssim\|u\|_{B^{N}_{1,1}}\nonumber,
\end{array}
$$
where we used the fact that for $q\geq0$, the homogeneous space
($\dot{\Delta}_{q}$) and the inhomogeneous space ($\Delta_{q}$)
share the same definition.

\textbf{Lemma 2.3}\quad Let $s>0,\ p,\ q\in[1,\infty]$, then
$B^{s}_{p,q}\bigcap L^{\infty}$ is an algebra and the following
inequality holds true
$$
\|uv\|_{B^{s}_{p,q}}\lesssim\|u\|_{L^{\infty}}\|v\|_{B^{s}_{p,q}}+\|v\|_{L^{\infty}}\|u\|_{B^{s}_{p,q}}.
$$
The proof can be found in the reference [10], and we omit it here.

\textbf{Lemma 2.4}\quad Let $r>0$, then there exists a constant $C$
such that the following inequality holds
$$
\|u\cdot\nabla v\|_{r}\leq C\|u\|_{r}\|v\|_{B^{1}_{\infty,1}}.
$$
\textbf{Proof.} We use Bony's formula for paraproduct of two
functions
$$
\begin{array}{llll}
u\cdot\nabla
v&=\displaystyle{T_{u^{i}}\partial_{i}v+T_{\partial_{i}v}u^{i}+R(u^{i},\partial_{i}v)}\\
&=\displaystyle{\sum_{k}S_{k-1}u^{i}\Delta_{k}\partial_{i}v+\sum_{k}S_{k-1}\partial_{i}v\Delta_{k}u^{i}}+\displaystyle{\sum_{k}(\Delta_{k}u^{i}\sum_{j=-1}^{1}\Delta_{k+j}\partial_{i}v)}.
\end{array}
$$
where we use the denotation of Einstein's sum about $i$.

From Proposition 2.1, we have
$$
\begin{array}{rl}
\Delta_{q}(u\cdot\nabla v)&=\displaystyle{\sum_{|k-q|\leq
M}\Delta_{q}(S_{k-1}u^{i}\Delta_{k}\partial_{i}v)+\sum_{|k-q|\leq
M}\Delta_{q}(S_{k-1}\partial_{i}v\Delta_{k}u^{i})}\\
&+\displaystyle{\sum_{k\geq
q-M}\Delta_{q}(\Delta_{k}u^{i}\sum^{1}_{j=-1}\Delta_{k+j}\partial_{i}v)}\\
&=L_{1}+L_{2}+L_{3}.
\end{array}
$$
where $M$ is a finite integer.\\
Thanks to H\"{o}lder inequality and Proposition 2.2, we obtain
\begin{equation}
\begin{array}{rl}
\displaystyle{\sup_{q\geq-1}2^{qr}\|L_{1}\|_{L^{\infty}}}&\leq\displaystyle{\sup_{q\geq-1}2^{qr}\sum_{|k-q|\leq
M}\|\Delta_{q}u^{i}\|_{L^{\infty}}\|\Delta_{k}\partial_{i}v\|_{L^{\infty}}}\\
&\leq C\displaystyle{\sup_{q\geq-1}2^{qr}\|\Delta_{q}u^{i}\|_{L^{\infty}}\sum_{k\geq-1}2^{k}\|\Delta_{k}v\|_{L^{\infty}}}\\
&\leq C\displaystyle{\|u\|_{r}\|v\|_{B^{1}_{\infty,1}}},
\end{array}
\end{equation}
and similarly
\begin{equation}
\begin{array}{rl}
\displaystyle{\sup_{q\geq-1}2^{qr}\|L_{2}\|_{L^{\infty}}}&\leq\displaystyle{\sup_{q\geq-1}2^{qr}\sum_{|k-q|\leq
M}\|\Delta_{q}\partial_{i}v\|_{L^{\infty}}\|\Delta_{k}u^{i}\|_{L^{\infty}}}\\
&\leq C\displaystyle{\sup_{q\geq-1}2^{qr}\|\Delta_{q}u^{i}\|_{L^{\infty}}\sum_{q\geq-1}2^{q}\|\Delta_{q}v\|_{L^{\infty}}}\\
&\leq\displaystyle{C\|u\|_{r}\|v\|_{B^{1}_{\infty,1}}}.
\end{array}
\end{equation}
On the other hand, we have
\begin{equation}
\begin{array}{rl}
\displaystyle{\sup_{q\geq-1}2^{qr}\|L_{3}\|_{L^{\infty}}}&\leq\displaystyle{\sup_{q\geq-1}2^{qr}\sum_{k\geq
q-M}\|\Delta_{q}u^{i}\|_{L^{\infty}}\|\sum^{1}_{j=-1}\Delta_{k+j}\partial_{i}v\|_{L^{\infty}}}\\
&\leq
\displaystyle{C\sup_{q\geq-1}2^{qr}\|\Delta_{q}u^{i}\|_{L^{\infty}}\sum_{k\geq-1}2^{k}\|v\|_{L^{\infty}}}\\
&\leq C\|u\|_{r}\|v\|_{B^{1}_{\infty,1}}.
\end{array}
\end{equation}
Using (2.1)-(2.3), we get the proof of Lemma 2.4.

As usual, the Riesz operator is defined as
$$
R_{j}(u)=\frac{x_{j}}{|x|^{N+1}}*u={\mathcal
{F}}^{-1}(-\frac{i\xi_{j}}{|\xi|}\cdot\hat{u}),
$$
where $i, j$ satisfy $i^{2}=-1$ and $1\leq j\leq N$. Hence
$$
R_{j}R_{k}(u)=\sum_{1\leq j,k\leq
N}\frac{\partial_{i}\partial_{j}}{(-\Delta)}(u).
$$

The following lemma is about the boundedness of the Riesz operator
in the space of $C^{r}(r>0)$ (see [8]).

\textbf{Lemma 2.5}\quad Let $r>0$, if $u\in C^{r}$, then for the
Riesz operator, there exists a constant $C$, such that
$$
\|R_{j}R_{k}(u)\|_{r}\leq C\|u\|_{r} \quad  i.e.\quad
\|\nabla\Delta^{-1}\mbox{div} u\|_{r}\leq C\|u\|_{r}.
$$

\textbf{Proof.} Note that
$$
\begin{array}{rl}
\|\nabla\Delta^{-1}\mbox{div}
u\|_{r}&=\displaystyle{\sup_{q\geq-1}2^{rq}\|\nabla\Delta^{-1}\mbox{div}
\Delta_{q}u\|_{L^{\infty}}}\nonumber\\
&\leq\displaystyle{\sup_{q\geq0}2^{rq}\|\nabla\Delta^{-1}\mbox{div}
\Delta_{q}u\|_{L^{\infty}}}
+\displaystyle{2^{-r}\|\nabla\Delta^{-1}\mbox{div}
\Delta_{-1}u\|_{L^{\infty}}}\nonumber\\
&\leq\displaystyle{\sup_{q\in Z}2^{rq}\|\nabla\Delta^{-1}\mbox{div}
\dot{\Delta}_{q}u\|_{L^{\infty}}+2^{-r}\|\nabla\Delta^{-1}\mbox{div}
\Delta_{-1}u\|_{L^{\infty}}}\nonumber\\
&=I+II\nonumber.
\end{array}
$$
Using Lemma 2.2(3), we have
\begin{equation}
I\leq\|u\|_{\dot{C}^{r}}\leq\|u\|_{r}.
\end{equation}
Since
$$
\Delta_{-1}u=\Delta_{-1}\sum_{q\in
Z}\dot{\Delta}_{q}u=\sum_{q\leq-1}\Delta_{-1}\dot{\Delta}_{q}u,
$$
we obtain
\begin{equation}
\begin{array}{rl}
II=\|\Delta_{-1}u\|_{L^{\infty}}&\leq\displaystyle{\sum_{q\leq-1}\|\dot{\Delta}_{q}u\|_{L^{\infty}}
=\|u\|_{\dot{B}^{0}_{\infty,1}}\lesssim\|u\|_{\dot{B}^{N}_{1,1}}}\\
&\lesssim\displaystyle{\|u\|_{B^{N}_{1,1}}\lesssim\|u\|_{r}},
\end{array}
\end{equation}
where we have used the boundedness property of the singular integral
operator (i.e. Riesz operator) from $\dot{B}^{s}_{p,q}$ into itself
in [5].

 According to (2.4) and (2.5), we have
$$
\|\nabla\Delta^{-1}\mbox{div} u\|_{r}\leq \|u\|_{r},
$$
The proof of the  lemma is complete.

Our main results of this paper is stated as

\vspace{3mm}

\textbf{Theorem 2.1(local existence and uniqueness)}  Suppose that
the initial data satisfy $\theta_{0}, u_{0}\in C^{r} (r>1)$. Then,
there exists $T^{*}=T(\|\theta_{0}\|_{C^{r}},\|u_{0}\|_{C^{r}})>0$,
such that the system (1.3)-(1.4) has a unique solution $(u,\theta)$
satisfying $u\in L^{\infty}([0,T^{*}];C^{r})$ and $\theta\in
L^{\infty}([0,T^{*}];C^{r})$.

\textbf{Theorem 2.2 (blow-up criterion)}\quad For $r>1$, if we
assume that the solution satisfies
$$
\int_{0}^{T^{*}}\|\nabla u\|_{L^{\infty}}<\infty,
$$
then the solution can be extended after $t=T^{*}$. In other word, if
the solution blows up at $t=T^{*}$, then
$$
\int_{0}^{T^{*}}\|\nabla u\|_{L^{\infty}}ds=\infty,
$$
for any pair of solution $(\theta,u)$ in the $C^{r}$ space.

\vspace{0.2cm}

\renewcommand{\thefootnote}{\fnsymbol{footnote}}
\thispagestyle{empty} \setcounter{section}{3}
\setcounter{equation}{0}
 { \large\bf 3. A Priori
Estimates} \vspace{0.2cm}

\vspace{3mm}

 In this section, we will prove the existence part of
Theorem 2.1. To this end,  we first derive some a priori estimates.

\textbf{Lemma 3.1} Let $r>0$, $v$ be a divergence-free vector field
belonging to the space $L^{1}_{loc}((0,+\infty);\mbox{Lip}(R^{N}))$
and $f$ be a scalar solution to the following problem
\begin{equation}
 \left\{
 \begin{array}{llll}
\partial_{t}f+v\cdot\nabla f=g,
\\
f|_{t=0}=f_{0}.
\end{array}\right.\nonumber
 \end{equation}
If the initial data $f_{0}\in C^{r}$, then we have for all $t\in
(0,+\infty)$
$$
\|f(t)\|_{r}\leq\|f_{0}\|_{r}+\int_{0}^{t}\|g(s)\|_{r}ds+C(r)\int_{0}^{t}\|\nabla
v\|_{L^{\infty}}\|f(s)\|_{r}ds,
$$
where $C$ depends only on the dimension $N$ and $r$.

\textbf{Proof.} Taking operation $\Delta_{q}$ on both sides of the
above system (3.1), we get
\begin{equation}
 \left\{
 \begin{array}{llll}
\partial_{t}\Delta_{q}f+v\cdot\nabla\Delta_{q}f=\Delta_{q}g+[v\cdot\nabla,\Delta_{q}]f,
\\
\Delta_{q}f|_{t=0}=\Delta_{q}f_{0}.
\end{array}\right.\nonumber
 \end{equation}
It's easy to get that
$$
\Delta_{q}f(t)=\Delta_{q}f_{0}+\int_{0}^{t}(\Delta_{q}g+[v\cdot\nabla,\Delta_{q}]f)ds.
$$
According to Lemma 2.1, we obtain
$$
\|f(t)\|_{r}\leq\|f_{0}\|_{r}+\int_{0}^{t}\|g(s)\|_{r}ds+C(r)\int_{0}^{t}\|\nabla
v\|_{L^{\infty}}\|f(s)\|_{r}ds.
$$

The proof of the lemma is finished.

Based on Lemma 3.1, we have

\textbf{Lemma 3.2} Let $r>0$. Suppose that $u ,\theta$ are smooth
solutions of (1.3) with initial data $u_0, \theta_0 \in C^{r}$. Then
we have
\begin{equation}
\|\theta(t)\|_{r}\leq\|\theta_{0}\|_{r}+C(r)\int_{0}^{t}\|\nabla
u\|_{L^{\infty}}\|\theta(s)\|_{r}ds.
 \end{equation}

\textbf{Remark 3.1} According to Lemma 3.2 and the Gronwall's
inequality, we have
\begin{equation}
\|\theta\|_{r}\leq\|\theta_{0}\|_{r}\exp(C(r)\int_{0}^{t}\|\nabla
u\|_{L^{\infty}}ds).
\end{equation}

\textbf{Lemma 3.3} Let $r>0$. Suppose that $u ,\theta$ are smooth
solutions of (1.3) with initial data $u_0, \theta_0 \in C^{r}$. Then
we have
\begin{equation}
\|u\|_{r}\leq\|u_{0}\|_{r}+2C(r)\int_{0}^{t}\|u\|_{r}\|\nabla
u\|_{L^{\infty}}ds+(2+2^{-r})\int_{0}^{t}\|\theta\|_{r}ds.
\end{equation}
\textbf{Proof.} In view of Lemma 3.1, we get
$$
\|u(t)\|_{r}\leq\|u_{0}\|_{r}+\int_{0}^{t}\|\nabla\Pi\|_{r}ds+\int_{0}^{t}\|\theta
e_{N}\|_{r}+C(r)\int_{0}^{t}\|\nabla u\|_{L^{\infty}}\|u(s)\|_{r}ds,
$$
where
$$
\nabla\Pi=-\nabla\Delta^{-1}\mbox{div}(u\cdot\nabla
u)+\nabla\Delta^{-1}\partial_{N}\theta.
$$
Note that
$$
\nabla\Delta_{q}\Pi=-\Delta_{q}\nabla\Delta^{-1}\mbox{div}(u\cdot\nabla
u)+\Delta_{q}\nabla\Delta^{-1}\partial_{N}\theta.
$$
Then one has
$$
\|\nabla\Delta\Pi\|_{L^{\infty}}\leq\|\Delta_{q}\nabla\Delta^{-1}\mbox{div}(u\cdot\nabla
u)\|_{L^{\infty}}+\|\Delta_{q}\nabla\Delta^{-1}\partial_{N}\theta\|_{L^{\infty}},
$$
and
$$
\sup_{q\geq-1}2^{qr}\|\nabla\Delta\Pi\|_{L^{\infty}}\leq\sup_{q\geq-1}2^{qr}\|\Delta_{q}\nabla\Delta^{-1}\mbox{div}(u\cdot\nabla
u)\|_{L^{\infty}}+\sup_{q\geq-1}2^{qr}\|\Delta_{q}\nabla\Delta^{-1}\partial_{N}\theta\|_{L^{\infty}}.
$$
Using Lemma 2.5 , Lemma 2.3 and Lemma 2.2 (3), one has (see [8]):
$$
\displaystyle{\sup_{q\geq-1}2^{qr}\|\Delta_{q}\nabla\Delta^{-1}\mbox{div}(u\cdot\nabla
u)\|_{L^{\infty}}\leq C(r)\|u\|_{r}\|\nabla u\|_{L^{\infty}}}.
$$

Concerning the term
$\displaystyle{\sup_{q\geq-1}2^{qr}\|\Delta_{q}\nabla\Delta^{-1}\partial_{N}\theta\|_{L^{\infty}}}$,
we have
\begin{equation}
\begin{array}{rl}
\displaystyle{\sup_{q\geq-1}2^{qr}\|\Delta_{q}\nabla\Delta^{-1}\partial_{N}\theta\|_{L^{\infty}}}
&\leq\displaystyle{\sup_{q\geq0}2^{qr}\|\Delta_{q}\nabla\Delta^{-1}\partial_{N}\theta\|_{L^{\infty}}+
2^{-r}\|\Delta_{-1}\nabla\Delta^{-1}\partial_{N}\theta\|_{L^{\infty}}}\\
&=I+II.
\end{array}
\end{equation}
Direct estimates give
\begin{equation}
\begin{array}{llll}
I&=\displaystyle{\sup_{q\geq0}2^{qr}\|\Delta_{q}\nabla\Delta^{-1}\partial_{N}\theta\|_{L^{\infty}}
\leq\sup_{q\in
Z}2^{qr}\|\dot{\Delta}_{q}\nabla\Delta^{-1}\partial_{N}\theta\|_{L^{\infty}}}\\
&\leq\displaystyle{\|\theta\|_{\dot{C}^{r}}\leq\|\theta\|_{r}},
\end{array}
\end{equation}
where we used the embedding
$B^{r}_{\infty,\infty}\hookrightarrow\dot{B}^{r}_{\infty,\infty}
(r>1)$ in Lemma 2.2 (3). Here we could also get (3.7) by using
directly Lemma 2.5.

By the dyadic decomposition in the homogeneous space, we have
\begin{equation}
\begin{array}{rl}
II&=2^{-r}\|\Delta_{-1}\nabla\Delta^{-1}\partial_{N}\theta\|_{L^{\infty}}=
\displaystyle{2^{-r}\|\Delta_{-1}\sum_{j\in
Z}\dot{\Delta}_{j}\nabla\Delta^{-1}\partial_{N}\theta\|_{L^{\infty}}}\nonumber\\
&\leq\displaystyle{2^{-r}\sum_{j\in
Z}\|\dot{\Delta}_{j}\nabla\Delta^{-1}\partial_{N}\theta\|_{L^{\infty}}}\nonumber\\
&\leq2^{-r}\|\theta\|_{\dot{B}^{0}_{\infty,1}}\leq2^{-r}\|\theta\|_{\dot{B}^{N}_{1,1}},
\end{array}
\end{equation}
where we used the boundedness of the Riesz operator in any
homogeneous Besov spaces $\dot{B}^{s}_{p,r}$ in [5].

Then in view of the embedding inequalities
$\|\theta\|_{\dot{B}^{N}_{1,1}}\lesssim\|\theta\|_{B^{N}_{1,1}}\lesssim\|\theta\|_{r}(r>0)$,
 we have $II\lesssim2^{-r}\|\theta\|_{r}$.

Putting (3.7) and (3.8) into (3.6), one has
$$
\sup_{q\geq-1}2^{qr}\|\triangle_{q}\nabla\triangle^{-1}\partial_{N}\theta\|_{L^{\infty}}\leq(1+2^{-r})\|\theta\|_{r}.
$$

So about the pressure term, we have
$$
\|\nabla\Pi\|_{r}=\sup_{q\leq-1}2^{qr}\|\nabla\Delta_{q}\Pi\|_{L^{\infty}}\leq
C(r)(\|u\|_{r}\|\nabla u\|_{L^{\infty}}+\|\theta\|_{r}).
$$
Thus we get
$$
\begin{array}{l}
\|u\|_{r}\leq\displaystyle{\|u_{0}\|_{r}+2C(r)\int_{0}^{t}\|u\|_{r}\|\nabla
u\|_{L^{\infty}}ds+(2+2^{-r})\int_{0}^{t}\|\theta\|_{r}ds}.
\end{array}
$$
The proof of the lemma is finished.

\vspace{3mm}

\thispagestyle{empty} \setcounter{section}{4}
\setcounter{equation}{0}
 { \large\bf 4. Proof of Main Results} \vspace{0.2cm}

\vspace{3mm}

 In this section we will prove one of our main results,  Theorem 2.1.

\vspace{3mm}
\thispagestyle{empty} \setcounter{subsection}{4}
\setcounter{equation}{0}
 { \large \bf 4.1 The approximate solutions} \vspace{0.2cm}

\vspace{3mm}

  In order to establish the local in time existence of
solution, we construct the approximate sequences
$\{(\theta^{n},u^{n})\}_{n\in Z^{+}\cup\{0\}}$ as follows
\begin{equation}
 \left\{
 \begin{array}{llll}
\partial_{t}\theta_{n+1}+u_{n}\cdot\nabla\theta_{n+1}=0,
\\
\partial_{t}u_{n+1}+u_{n}\cdot\nabla u_{n+1}+\nabla\Pi_{n+1}=\theta_{n+1} e_{N},
 \\
{\rm div} u_{n}=0={\rm div}u_{n+1},
\end{array}\right.
 \end{equation}
with the initial data
\begin{equation}
 \left\{
 \begin{array}{llll}
\theta_{1}=S_{2}\theta_{0},\\
\theta_{n+1}|_{t=0}=S_{n+2}\theta_{0},\\
u_{1}=S_{2}u_{0},\\
u_{n+1}|_{t=0}=S_{n+2}u_{0}.
 \end{array}\right.
\end{equation}
Taking  $\Delta_{q}$ on both sides of (4.1)$_{1}$ and (4.2)$_{2}$,
we have
\begin{equation}
 \left\{
 \begin{array}{llll}
\partial_{t}\Delta_{q}\theta_{n+1}+u_{n}\cdot\nabla\Delta_{q}\theta_{n+1}=[u_{n}\cdot\nabla,\Delta_{q}]\theta_{n+1},\\
\Delta_{q}\theta_{n+1}|_{t=0}=\Delta_{q}S_{n+2}\theta_{0}.\nonumber
\end{array}\right.
 \end{equation}
By Lemma 3.1, we get
$$
\begin{array}{l}
\displaystyle{\|\theta_{n+1}\|_{r}\leq\|\theta_{0}\|_{r}+C(r)\int_{0}^{t}\|\theta_{n+1}\|_{r}\|\nabla
u_{n}\|_{L^{\infty}}ds.}
\end{array}
$$
Using the  Gronwall's inequality, we have
\begin{equation}
\|\theta_{n+1}(t)\|_{r}\leq\|\theta_{0}\|_{r}\exp(C(r)\int_{0}^{t}\|\nabla
u_{n}\|_{L^{\infty}}ds).
\end{equation}
Taking  $\Delta_{q}$ on both sides of (4.1)$_{2}$ and (4.2)$_{4}$,
we get
$$
 \left\{
 \begin{array}{llll}
\partial_{t}\Delta_{q}u_{n+1}+u_{n}\cdot\nabla\Delta_{q}u_{n+1}=[u_{n}\cdot\nabla,\Delta_{q}]u_{n+1}-\nabla\Delta_{q}\Pi_{n+1}+\Delta_{q}\theta_{n}e_{N}\nonumber,\\
\Delta_{q}u_{n+1}|_{t=0}=\Delta_{q}S_{n+2}u_{0}.\nonumber
\end{array}\right.
$$
Similar to Lemma 3.1, we get
$$
\Delta_{q}u_{n+1}(t)=\Delta_{q}u_{n+1}(t=0)+\int_{0}^{t}([u_{n}\cdot\nabla,\Delta_{q}]u_{n+1}-\nabla\Delta_{q}\Pi_{n+1}+\Delta_{q}\theta_{n}e_{N})ds,
$$
here the pressure term
$$
\nabla\Pi_{n+1}=-\nabla\Delta^{-1}\mbox{div}(u_{n}\cdot\nabla
u_{n+1})+\nabla\Delta^{-1}\partial_{N}\theta_{n}.
$$
So we have
\begin{equation}
\begin{array}{rl}
\|u_{n+1}(t)\|_{r}&\leq\|u_{n+1}(t=0)\|_{r}+\displaystyle{\int_{0}^{t}(\sup_{q\geq-1}2^{qs}
\|[u_{n}\cdot\nabla,\Delta_{q}]u_{n+1}\|_{L^{\infty}}\nonumber}\\
&+\displaystyle{\sup_{q\geq-1}2^{qs}\|\nabla\Delta_{q}\Pi_{n+1}\|_{L^{\infty}}
+\sup_{q\geq-1}2^{qs}\|\Delta_{q}\theta_{n+1}e_{N}\|_{L^{\infty}})ds}.
\end{array}
\end{equation}
In view of Lemma 2.1, the following estimates hold
$$
\sup_{q\geq-1}2^{qs}\|[u_{n}\cdot\nabla,\Delta_{q}]u_{n+1}\|_{L^{\infty}}\leq
C(r)\|\nabla u_{n}\|_{L^{\infty}}\|u_{n+1}\|_{r},
$$
and
$$
\sup_{q\geq-1}2^{qs}\|\Delta_{q}\theta_{n+1}e_{N}\|_{L^{\infty}}\leq\|\theta_{n+1}\|_{r}.
$$
Concerning  the pressure term, we have
$$
\begin{array}{rl}
\|\nabla\Pi_{n+1}\|_{r}&=\displaystyle{\sup_{q\geq-1}2^{qs}\|\nabla\Delta_{q}\Pi_{n+1}\|_{L^{\infty}}
\leq\sup_{q\geq-1}2^{qs}\|\Delta_{q}\nabla\Delta^{-1}\mbox{div}(u_{n}\cdot\nabla
u_{n+1})\|_{L^{\infty}}}\nonumber\\
&+\displaystyle{\sup_{q\geq-1}2^{qs}\|\Delta_{q}\nabla\Delta^{-1}\partial_{N}\theta_{n+1}\|_{L^{\infty}}=J_{1}+J_{2}}\nonumber.
\end{array}
$$
Similar as in [8] Proposition 2.5.1, we obtain
\begin{equation}
\begin{array}{rl}
J_{1}&=\displaystyle{\sup_{q\geq-1}2^{qs}\|\Delta_{q}\nabla\Delta^{-1}\mbox{div}(u_{n}\cdot\nabla
u_{n+1})\|_{L^{\infty}}}\nonumber\\
&\leq C(r)(\|\nabla u_{n}\|_{L^{\infty}}\|u_{n+1}\|_{r}+\|\nabla
u_{n+1}\|_{L^{\infty}}\|u_{n}\|_{r}),
\end{array}
\end{equation}
and
\begin{equation}
\begin{array}{rl}
J_{2}&\leq\displaystyle{\sup_{q\in
Z}\|\dot{\Delta}_{q}\nabla\Delta^{-1}\partial_{N}\theta_{n+1}\|+2^{-r}\|\Delta_{-1}\nabla\Delta^{-1}\partial_{N}\theta_{n+1}\|}\nonumber\\
&\leq\|\theta_{n+1}\|_{\dot{C}^{r}}+\|\theta_{n+1}\|_{\dot{B}^{0}_{\infty,1}}\lesssim\|\theta_{n+1}\|_{r}.
\end{array}
\end{equation}
Putting (4.6) and (4.7) into (4.5), we have
$$
\begin{array}{rl}
\|u_{n+1}(t)\|_{r}&\leq\|u_{n+1}(0)\|_{r}\nonumber\\
&+\displaystyle{\int_{0}^{t}(2C(r)\|\nabla
u_{n}\|_{L^{\infty}}\|u_{n+1}\|_{r}+C(r)\|\nabla
u_{n+1}\|_{L^{\infty}}\|u_{n}\|_{r}+2\|\theta_{n}\|_{r})ds}\nonumber\\
&\leq\displaystyle{\|u_{n+1}(0)\|_{r}+3C(r)\int_{0}^{t}\|u_{n}\|_{r}\|u_{n+1}\|_{r}ds+2\int_{0}^{t}\|\theta_{n}\|_{r}ds}\nonumber
\end{array}
$$
Applying the Gronwall's inequality yields
$$
\|u_{n+1}\|_{r}\leq[\|S_{n+2}u_{0}\|_{r}+2\int_{0}^{t}\|\theta_{n+1}\|_{r}ds]\exp(\int_{0}^{t}3C(r)\|u_{n}\|_{r}ds),
$$
i.e. \begin{equation}
\|u_{n+1}\|_{r}\leq[\|u_{0}\|_{r}+2\int_{0}^{t}\|\theta_{n+1}\|_{r}ds]\exp(\int_{0}^{t}3C(r)\|u_{n}\|_{r}ds).
\end{equation}

Now we define $a_{0}=\|\widetilde{h}\|_{L^{1}}$(the function
$\widetilde{h}$ is given in the notation).  We will establish that,
for all initial data $\theta_{0},u_{0}$,
\begin{equation}
\|\theta_{n}\|_{C([0,T_{1}];C^{r})}\leq Pa_{0}\|\theta_{0}\|_{r}
\quad \|u_{n}\|_{C([0,T_{1}];C^{r})}\leq Qa_{0}\|u_{0}\|_{r}
\end{equation}
 with some constants $P$ and $Q$ which may be large enough ( for example, $P=Q=32$), and
$T_{1}>0$ which will be determined later.

Using  (4.4) and (4.8), we have
 $$
\|\theta_{1}\|_{r}\leq\|\theta_{0}\|_{r}\exp(C(r)\int_{0}^{t}\|u_{0}\|_{r}ds)\leq
Pa_{0}\|\theta_{0}\|_{r}.
 $$
Let
$$
\displaystyle{T_{1}^{(1)}=\frac{1}{C(r)\|u_{0}\|_{r}}\ln(Pa_{0})},
$$
Then we obtain that
$\displaystyle{\|\theta_{1}\|_{C([0,T_{1}^{(1)}],C^{r})}\leq
Pa_{0}\|\theta_{0}\|_{r}}$, when $t<T_{1}^{(1)}$.

Then we let
$$
\begin{array}{rl}
\|u_{1}\|_{r}&\leq\displaystyle{\|u_{0}\|_{r}[1+2tPa_{0}\frac{\|\theta_{0}\|_{r}}{\|u_{0}\|_{r}}]
\exp(3tC(r)\|u_{0}\|_{r})}\nonumber\\
&\leq\displaystyle{\|u_{0}\|_{r}\exp(3tC(r)\|u_{0}\|_{r}+2tPa_{0}\frac{\|\theta_{0}\|_{r}}{\|u_{0}\|_{r}})
\leq Qa_{0}\|u_{0}\|_{r}}\nonumber,
\end{array}
$$
Let
$$
\displaystyle{T_{1}^{(2)}=\frac{\ln(Qa_{0})}{3C(r)\|u_{0}\|_{r}+2Pa_{0}\frac{\|\theta_{0}\|_{r}}{\|u_{0}\|_{r}}}}.
$$
Then we get $\displaystyle{\|u_{1}\|_{C([0,T_{1}^{(2)}],C^{r})}\leq
Qa_{0}\|u_{0}\|_{r}}$, when $t\leq T_{1}^{(2)}$.

Now we apply the induction of $n$ to obtain (4.9). We assume the
 estimates (4.9) to be true for every $j\leq n$.

Firstly, we have
$$
\|\theta_{n+1}\|_{r}\leq\|\theta_{0}\|_{r}\exp(C(r)\int_{0}^{t}\|u_{n}\|_{r}ds)\leq\|\theta_{0}\|_{r}\exp(tC(r)\|u_{n}\|_{r})\leq
Pa_{0}\|\theta_{0}\|_{r}.
$$
Then we set
$$
\displaystyle{T_{1}^{(3)}=\frac{\ln(Pa_{0})}{C(r)Qa_{0}\|u_{0}\|_{r}}}.
$$
If $t\leq T_{1}^{(3)}$, the inequality $
\|\theta_{n+1}\|_{C([0,T_{1}^{(3)}],C^{r})}\leq
Pa_{0}\|\theta_{0}\|_{r} $ is obtained.

Then for the term $u_{n+1}$, we let
$$
\begin{array}{rl}
\|u_{n+1}\|_{r}&\leq\displaystyle{[\|u_{0}\|_{r}+2\int_{0}^{t}\|\theta_{n+1}\|_{r}ds]\exp(\int_{0}^{t}3C(r)\|u_{n}\|_{r}ds)}\nonumber\\
&\leq\displaystyle{[\|u_{0}\|_{r}+2tPa_{0}\|\theta_{0}\|_{r}]\exp(3tC(r)Qa_{0}\|u_{0}\|_{r})}\nonumber\\
&\leq\displaystyle{\|u_{0}\|_{r}\exp\|u_{0}\|_{r}\exp(3tC(r)Qa_{0}\|u_{0}\|_{r}+2tPa_{0}\frac{\|\theta_{0}\|_{r}}{\|u_{0}\|_{r}})}
\leq Qa_{0}\|u_{0}\|_{r}.\nonumber
\end{array}
$$
Set
$$
\displaystyle{T_{1}^{(4)}=\frac{\ln(Qa_{0})}{3C(r)Qa_{0}\|u_{0}\|_{r}+2Pa_{0}\frac{\|\theta_{0}\|_{r}}{\|u_{0}\|_{r}}}}.
$$
If $t\leq T_{1}^{(4)}$, there holds
$\displaystyle{\|u_{n+1}\|_{C([0,T_{1}^{(4)}],C^{r})}\leq
Qa_{0}\|u_{0}\|_{r}}$.

Thus the sequences $\{(\theta_{n},u_{n})\}_{n\in Z^{+}\cup\{0\}}$
are constructed, and  they are also bounded in the space
$C([0,T_{1}],C^{r})$, where $T_{1}=\min\{T_{1}^{i}\}_{i=1}^{4}$.

\vspace{3mm}

\thispagestyle{empty} \setcounter{subsection}{4}
 { \large \bf 4.2 The Cauchy sequences} \vspace{0.2cm}

\vspace{3mm}

  Now we prove that there exists a $T_{2}>0$, such that the
sequences $\displaystyle{(\theta_{n})_{n\in Z^{+}\cup\{0\} }}$,

\noindent$\displaystyle{(u_{n})_{n\in Z^{+}\cup\{0\} }}$ are the
Cauchy sequences in the space $C([0,T_{2}];C^{r-1})(r>1)$. To do so,
we  estimate the quantities
$\|\theta_{n+1}(t)-\theta_{n}(t)\|_{r-1}$ and
$\|u_{n+1}(t)-u_{n}(t)\|_{r-1}$. For conciseness,  we set
$\bar{u}_{n+1}(t)=u_{n+1}(t)-u_{n}(t)$ and
$\bar{\theta}_{n+1}(t)=\theta_{n+1}(t)-\theta_{n}(t)$.

By constructions of the approximate solutions, we know that
\begin{equation}
 \left\{
 \begin{array}{llll}
\partial_{t}\theta_{n+1}+u_{n}\cdot\nabla\theta_{n+1}=0,
\\
\theta_{n+1}|_{t=0}=S_{n+2}\theta_{0}, \nonumber
\end{array}\right.
 \end{equation}
and
\begin{equation}
 \left\{
 \begin{array}{llll}
\partial_{t}\theta_{n}+u_{n-1}\cdot\nabla\theta_{n}=0,
\\
\theta_{n}|_{t=0}=S_{n+1}\theta_{0}. \nonumber
\end{array}\right.
 \end{equation}
Subtracting (4.11) from (4.10), we get
\begin{equation}
 \left\{
 \begin{array}{llll}
\partial_{t}(\theta_{n+1}-\theta_{n})+u_{n}\cdot\nabla(\theta_{n+1}-\theta_{n})=-(u_{n}-u_{n-1}\cdot\nabla\theta_{n}),
\\
(\theta_{n+1}-\theta_{n})|_{t=0}=S_{n+2}\theta_{0}-S_{n+1}\theta_{0}=\Delta_{n+1}\theta_{0}.
\nonumber
\end{array}\right.
 \end{equation}
In view of Lemma 3.2, we obtain
$$
\begin{array}{rl}
\|\bar{\theta}_{n+1}\|_{r-1}&\leq\displaystyle{\|\Delta_{n+1}\theta_{0}\|_{r-1}+\int_{0}^{t}
\|\bar{u}_{n}\cdot\nabla\theta_{n}\|_{r-1}ds}\nonumber\\
&+\displaystyle{C(r)\int_{0}^{t}\|\nabla
u_{n}\|_{L^{\infty}}\|\bar{\theta}_{n+1}\|_{r-1}ds}\nonumber.
\end{array}
$$
Thanks to Lemma 2.3, we have
$$
\begin{array}{rl}
\|\bar{u}_{n}\cdot\nabla\theta_{n}\|_{r-1}&\leq
\displaystyle{C(\|\bar{u}_{n}\|_{L^{\infty}}\|\nabla\theta_{n}\|_{r-1}+\|\nabla\theta_{n}\|_{L^{\infty}}\|\bar{u}_{n}\|_{r-1})}\nonumber\\
&\leq\displaystyle{C(\|\bar{u}_{n}\|_{r-1}\|\theta_{n}\|_{r}+\|\theta_{n}\|_{r}\|\bar{u}_{n}\|_{r-1})}\nonumber\\
&\leq
\displaystyle{2C\|\theta_{n}\|_{r}\|\bar{u}_{n}\|_{r-1}}\nonumber,
\end{array}
$$
where we used the embedding inequalities: $\|\nabla
u\|_{B^{s-1}_{p,q}}\lesssim\|u\|_{B^{s}_{p,q}}$ and

\noindent$\|\nabla
u\|_{L^{\infty}}\lesssim\|u\|_{\dot{B}^{1}_{\infty,1}}\lesssim\|u\|_{r},(r>1)$.

Therefore we get the estimate
$$
\begin{array}{rl}
\|\bar{\theta}_{n+1}\|_{r-1}&\leq\displaystyle{\|\Delta_{n+1}\theta_{0}\|_{r-1}+\int_{0}^{t}
2C\|\theta_{n}\|_{r}\|\bar{u}_{n}\|_{r-1}ds}\nonumber\\
&+\displaystyle{C(r)\int_{0}^{t}\|u_{n}\|_{r}\|\bar{\theta}_{n+1}\|_{r-1}ds}.
\end{array}
$$
Since
$$
\begin{array}{rl}
\|\Delta_{n+1}\theta_{0}\|_{r-1}&=\displaystyle{\sup_{q\in
Z}2^{q(r-1)}\|\dot{\Delta}_{q}\Delta_{n+1}\theta_{0}\|_{L^{\infty}}=\sup_{|q-n|\leq3}2^{q(r-1)}\|\dot{\Delta}_{q}\theta_{0}\|_{L^{\infty}}}\nonumber\\
&=\displaystyle{\sup_{|q-n|\leq3}2^{qr}\|\dot{\Delta}_{q}\theta_{0}\|_{L^{\infty}}2^{-q}\lesssim2^{-n}\|\theta_{0}\|_{r}}\nonumber,
\end{array}
$$
we obtain
\begin{equation}
\begin{array}{rl}
\|\bar{\theta}_{n+1}\|_{r-1}&\leq\displaystyle{2^{-n}\|\theta_{0}\|_{r}+\int_{0}^{t}
2C\|\theta_{n}\|_{r}\|\bar{u}_{n}\|_{r-1}ds}\\
&+\displaystyle{C(r)\int_{0}^{t}\|u_{n}\|_{r}\|\bar{\theta}_{n+1}\|_{r-1}ds}.
\end{array}
\end{equation}
Using the Gronwall's inequality, we have
\begin{equation}
\begin{array}{rl}
\|\bar{\theta}_{n+1}\|_{r-1}&\leq\displaystyle{2^{-n}\|\theta_{0}\|_{r}\exp(C(r)\int_{0}^{t}\|u_{n}\|_{r}ds)}\nonumber\\
&+\displaystyle{2C\int_{0}^{t}\|\theta_{n}\|_{r}\|\bar{u}_{n}\|_{r-1}ds\exp(C(r)\int_{0}^{t}\|u_{n}\|_{r}ds)}.
\end{array}
\end{equation}

Now we estimate $\|\bar{u}_{n+1}\|_{r-1}$. Noting that
\begin{equation}
 \left\{
 \begin{array}{llll}
\partial_{t}u_{n+1}+u_{n}\cdot\nabla u_{n+1}=-\nabla\Pi_{n+1}+\theta_{n+1}e_{N},
\\
u_{n+1}|_{t=0}=S_{n+2}u_{0}, \nonumber\\
{\rm div}u_{n+1}=0,\nonumber
\end{array}\right.
 \end{equation}
and that
\begin{equation}
 \left\{
 \begin{array}{llll}
\partial_{t}u_{n}+u_{n-1}\cdot\nabla u_{n}=-\nabla\Pi_{n}+\theta_{n}e_{N},
\\
u_{n}|_{t=0}=S_{n+1}u_{0}, \nonumber\\
{\rm div}u_{n}=0,\nonumber
\end{array}\right.
 \end{equation}
subtracting (4.16) from (4.15), we get
\begin{equation}
 \left\{
 \begin{array}{llll}
\partial_{t}\bar{u}_{n+1}+u_{n}\cdot\nabla\bar{u}_{n+1}=-\bar{u}_{n}\cdot\nabla u_{n}-\nabla\Pi_{n+1}+\nabla\Pi_{n}
+\bar{\theta}_{n+1}e_{N},
\\
\bar{u}_{n+1}|_{t=0}=S_{n+2}u_{0}-S_{n+1}u_{0}=\Delta_{n+1}u_{0}.
\nonumber
\end{array}\right.
 \end{equation}
In view of Lemma 3.2, we obtain
$$
\begin{array}{rl}
\|\bar{u}_{n+1}\|_{r-1}&\leq\displaystyle{\|\Delta_{n+1}u_{0}\|_{r-1}+\int_{0}^{t}
\|\nabla\Pi_{n+1}-\nabla\Pi_{n}\|_{r-1}ds+\int_{0}^{t}\|\bar{\theta}_{n+1}\|_{r-1}ds}\nonumber\\
&+\displaystyle{C(r)\int_{0}^{t}\|\nabla
u_{n}\|_{L^{\infty}}\|\bar{u}_{n+1}\|_{r-1}ds}\nonumber,
\end{array}
$$
where
$$
\nabla\Pi_{n+1}=-\nabla\Delta^{-1}{\rm div}(u_{n}\cdot\nabla
u_{n+1})+\nabla\Delta^{-1}\partial_{N}\theta_{n+1},
$$
$$
\nabla\Pi_{n}=-\nabla\Delta^{-1}{\rm div}(u_{n-1}\cdot\nabla
u_{n})+\nabla\Delta^{-1}\partial_{N}\theta_{n}.
$$
It follows that
$$
\begin{array}{rl}
\nabla\Pi_{n+1}-\nabla\Pi_{n}&=-\nabla\Delta^{-1}{\rm
div}(u_{n}\cdot\nabla
u_{n+1})+\nabla\Delta^{-1}\partial_{N}\theta_{n+1}\nonumber\\
&+\nabla\Delta^{-1}{\rm div}(u_{n-1}\cdot\nabla
u_{n})-\nabla\Delta^{-1}\partial_{N}\theta_{n}\nonumber\\
&=-\nabla\Delta^{-1}{\rm
div}(u_{n}\cdot\nabla\bar{u}_{n+1})-\nabla\Delta^{-1}{\rm
div}(\bar{u}_{n}\cdot\nabla
u_{n})+\nabla\Delta^{-1}\partial_{N}\bar{\theta}_{n+1}\nonumber
\end{array}
$$

When $1<r<2$, we use an estimate in [8], which says that, if
$r\in(-1,1)$, then
\begin{equation}
\|\pi(v,w)\|_{r}\leq
C(\frac{1}{1+r}+\frac{1}{1-r})\min\{\|v\|_{Lip}\|w\|_{r},\|v\|_{r}\|w\|_{Lip}\},
\end{equation}
where  $\pi$ is viewed as the term
$-\nabla\Delta^{-1}\mbox{div}(u\cdot\nabla u)$, i.e.
$\pi=\nabla\Pi-\nabla\Delta^{-1}\partial_{N}\theta$.

It follows from (4.18) that
\begin{equation}
\|\nabla\Pi_{n+1}-\nabla\Pi_{n}\|_{r-1} \leq
\displaystyle{C(\|\bar{u}_{n+1}\|_{r-1}\|u_{n}\|_{r}+\|\bar{u}_{n}\|_{r-1}\|u_{n}\|_{r})+\|\bar{\theta}_{n+1}\|_{r-1}}.
\end{equation}

When $r>2$, using the inequality $\|\nabla
u\|_{L^{\infty}}\leq\|u\|_{r-1}$, and using Lemma 2.3 and Lemma 2.5,
we get
$$
\begin{array}{rl}
\|\nabla\Pi_{n+1}-\nabla\Pi_{n}\|_{r-1}&\leq\displaystyle{\|u_{n}\cdot\nabla
\bar{u}_{n+1}\|_{r-1}+\|\bar{u}_{n}\cdot\nabla
u_{n}\|_{r-1}+\|\bar{\theta}_{n+1}\|_{r-1}}\nonumber\\
&\leq
\displaystyle{C(\| u_{n}\|_{L^{\infty}}\|\nabla\bar{u}_{n+1}\|_{r-1}+\|\nabla\bar{u}_{n+1}\|_{L^{\infty}}\|u_{n}\|_{r-1}}\nonumber\\
&+\displaystyle{\|\bar{u}_{n}\|_{L^{\infty}}\|\nabla
u_{n}\|_{r-1}+\|\nabla
u_{n}\|_{L^{\infty}}\|\bar{u}_{n}\|_{r-1})+\|\bar{\theta}_{n+1}\|_{r-1}}\nonumber.
\end{array}
$$
we can also obtain (4.19).

For $r=2$, using Lemma 2.4 and Lemma 2.3, we have
$$
\begin{array}{rl}
\|\nabla\Pi_{n+1}-\nabla\Pi_{n}\|_{1}\leq&\displaystyle{\|u_{n}\cdot\nabla\bar{u}_{n+1}\|_{1}+\|\bar{u}_{n}\cdot\nabla
u_{n}\|_{1}+\|\bar{\theta}_{n+1}\|_{1}}\\
&\leq
\displaystyle{C\|u_{n}\|_{1}\|\bar{u}_{n+1}\|_{B^{1}_{\infty,1}}+\|\bar{\theta}_{n+1}\|_{1}}\\
&+\displaystyle{C(\|\bar{u}_{n}\|_{L^{\infty}}\|\nabla u_{n}\|_{1}+\|\bar{u}_{n}\|_{1}\|\nabla u_{n}\|_{L^{\infty}})}\\
&\leq\displaystyle{C(\|\bar{u}_{n+1}\|_{1}\|u_{n}\|_{2}+\|\bar{u}_{n}\|_{1}\|u_{n}\|_{2})+\|\bar{\theta}_{n+1}\|_{1}}.
\end{array}
$$
(4.19) is proved  for $r=2$. And we have proved that  for all $r>1$,
(4.19) holds true.

According to the inequality about the pressure term and some
embedding inequalities, we have

$$
\begin{array}{rl}
\displaystyle{\|\bar{u}_{n+1}\|_{r-1}}&\leq\displaystyle{\|\Delta_{n+1}u_{0}\|_{r-1}+C(r)\int_{0}^{t}\|u_{n}\|_{r}\|\bar{u}_{n}\|_{r-1}ds}\nonumber\\
&+\displaystyle{2C(r)\int_{0}^{t}\|u_{n}\|_{r}\|\bar{u}_{n+1}\|_{r-1}ds+2\int_{0}^{t}\|\bar{\theta}_{n+1}\|_{r-1}ds}.
\end{array}
$$
Since $\|\Delta_{n+1}u_{0}\|_{r-1}\leq C2^{-n}\|u_{0}\|_{r}$, we
deduce that
$$
\begin{array}{rl}
\|\bar{u}_{n+1}\|_{r-1}&\leq
\displaystyle{C2^{-n}(\|u_{0}\|_{r}+\|\theta_{0}\|_{r}\exp(C(r)\int_{0}^{t}\|u_{n}\|_{r}ds))}\nonumber\\
&+\displaystyle{C(r)\int_{0}^{t}\|u_{n}\|_{r}\|\bar{u}_{n}\|_{r-1}ds
+2C(r)\int_{0}^{t}\|u_{n}\|_{r}\|\bar{u}_{n+1}\|_{r-1}ds}\nonumber\\
&+\displaystyle{2C\int_{0}^{t}\int_{0}^{\tau}\|\theta_{n}\|_{r}\|\bar{u}_{n}\|_{r-1}ds\exp(C(r)\int_{0}^{\tau}\|u_{n}\|_{r}ds)d\tau}\nonumber.
\end{array}
$$
Then we use the Gronwall's inequality to give that
$$
\begin{array}{rl}
&\displaystyle{\|\bar{u}_{n+1}\|_{C([0,T_{2}];C^{r-1}})}\\ &\leq
\displaystyle{C2^{-n}(\|u_{0}\|_{r}+\|\theta_{0}\|_{r}\exp(C(r)\int_{0}^{t}\|u_{n}\|_{r}ds))\exp(2C(r)\int_{0}^{t}\|u_{n}\|_{r}ds)}\nonumber\\
&+\displaystyle{C(r)\int_{0}^{t}\|u_{n}\|_{r}\|\bar{u}_{n}\|_{r-1}ds\exp(2C(r)\int_{0}^{t}\|u_{n}\|_{r}ds)}\nonumber\\
&+\displaystyle{2C\int_{0}^{t}\int_{0}^{\tau}\|\theta_{n}\|_{r}\|\bar{u}_{n}\|_{r-1}ds\exp(C(r)\int_{0}^{t}\|u_{n}\|_{r}ds)d\tau\exp(2C(r)\int_{0}^{t}\|u_{n}\|_{r}ds)}\nonumber\\
&\leq \displaystyle{C2^{-n}(\|u_{0}\|_{r}+\|\theta_{0}\|_{r}\exp(C(r)tQa_{0}\|u_{0}\|_{r}))\exp(2C(r)tQa_{0}\|u_{0}\|_{r})}\\
&+\displaystyle{\|\bar{u}_{n}\|_{C([0,T_{2}];C^{r-1})}(C(r)Qa_{0}\|u_{0}\|_{r}\exp(2C(r)tQa_{0}\|u_{0}\|_{r})}\\
&+\displaystyle{\frac{2CPa_{0}\|\theta_{0}\|_{r}}{C(r)Qa_{0}\|u_{0}\|_{r}}\exp(C(r)tQa_{0}\|u_{0}\|_{r})t\exp(2C(r)tQa_{0}\|u_{0}\|_{r})}\\
&+\displaystyle{\frac{2CPa_{0}\|\theta_{0}\|_{r}}{(C(r)Qa_{0}\|u_{0}\|_{r})^{2}}[\exp(C(r)tQa_{0}\|u_{0}\|_{r})-1])\exp(2C(r)tQa_{0}\|u_{0}\|_{r})}\\
&=\displaystyle{\sum_{i=1}^{4}I_{i}}.
\end{array}
$$

Now we deal with the four terms $I_i(i=1,2,3,4)$  one by one.

Concerning $I_{1}$,  we choose
$$
T_{2}^{(1)}=\min\{\frac{1}{2C(r)Qa_{0}\|u_{0}\|_{r}}\ln(\frac{S}{\|u_{0}\|_{r}}),\frac{1}{3C(r)Qa_{0}\|u_{0}\|_{r}}\ln(\frac{S}{\|\theta_{0}\|_{r}})\},
$$ to obtain
$$
 \|u_{0}\|_{r}\exp(2tC(r)Qa_{0}\|u_{0}\|_{r})\leq S
\quad\mbox{and}\quad
 \|\theta_{0}\|_{r}\exp(3tC(r)Qa_{0}\|u_{0}\|_{r})\leq S
$$
where $S$ is a real constant large enough.\\

Then we get $I_{1}\leq CS2^{-n+1}$.

Concerning $I_{2}$, we choose
$$
T_{2}^{(2)}=\frac{1}{3C(r)Qa_{0}\|u_{0}\|_{r}}\ln(\frac{Pa_{0}}{5C(r)Q\|u_{0}\|_{r}}),
$$
to get
$$
C(r)Qa_{0}\|u_{0}\|_{r}\exp(3C(r)tQa_{0}\|u_{0}\|_{r})\leq\frac{1}{5}Pa_{0}^{2}.
$$
Then we get $I_{2}\leq(1/5)\|\bar{u}\|_{C([0,T_{2}];C^{r-1})}$.

Concerning $I_{3}$, we choose
$$
\frac{2CPa_{0}\|\theta_{0}\|_{r}}{C(r)Qa_{0}\|u_{0}\|_{r}}\exp(3C(r)tQa_{0}\|u_{0}\|_{r})t\leq
\frac{1}{5}Qa_{0}^{2}.
$$
Then there exists  a constant $T_{2}^{(3)}$ such that when $t\leq
T_{2}^{(3)}$ the inequality
$I_{3}\leq(1/5)\|\bar{u}\|_{C([0,T_{2}];C^{r-1})}$ holds.

Concerning $I_{4}$, we set
$$
T_{2}^{(4)}=\frac{1}{C(r)Qa_{0}\|u_{0}\|_{r}}\ln(1+\frac{Pa_{0}(C(r)Qa_{0}\|u_{0}\|_{r})^{2}}{5CPa_{0}\|\theta_{0}\|_{r}}),
$$
such that
$$
\frac{2CPa_{0}\|\theta_{0}\|_{r}}{(C(r)Qa_{0}\|u_{0}\|_{r})^{2}}[\exp(C(r)tQa_{0}\|u_{0}\|_{r})-1]\leq
\frac{1}{5}Pa_{0}.
$$
Then we   get $I_{4}\leq(1/5)\|\bar{u}\|_{C([0,T_{2}];C^{r-1})}$.

Choosing $\displaystyle{T_{2}=\min\{T_{2}^{(i)}\}_{i=1}^{4}}$, we
obtain
$$
\|u_{n+1}-u_{n}\|_{C([0,T_{2}];C^{r-1})}\lesssim
2^{-n}+\frac{3}{5}\|u_{n}-u_{n-1}\|_{C([0,T_{2}];C^{r-1})}.
$$

Therefore, for $r>1$, the sequence $(u_{n})_{n\in Z^{+}\cup\{0\}}$
is a Cauchy sequence in $C([0,T_{2}];C^{r-1})$. Furthermore using
(4.14), we obtain that the sequence $(\theta_{n})_{n\in
Z^{+}\cup\{0\}}$ is a Cauchy sequence in $C([0,T_{2}];C^{r-1})$.

Let $T^{*}=\min\{T_{1},T_{2}\}$, and denote the limit of sequences
$(\theta_{n})_{n\in Z^{+}\cup\{0\}}, (u_{n})_{n\in Z^{+}\cup\{0\}}$
by $\theta(t,x), u(t,x)$ respectively, we obtain that
$\theta(t,x)\in L([0,T^{*}];C^{r})$, $u(t,x)\in L([0,T^{*}];C^{r})$,
where $r>1$, are solutions of (1.3)-(1.4).

\vspace{3mm}
\thispagestyle{empty} \setcounter{subsection}{4}
 { \large \bf 4.3 Uniqueness} \vspace{0.2cm}
\vspace{3mm}

 Suppose that $(\theta_{i},u_{i})\in
L^{\infty}([0,T];C^{r}(R^{N};R^{N}))$ $(i=1,2)$ are two solutions of
the system (1.3)-(1.4).

We set
$\theta=\theta_{1}-\theta_{2},u=u_{1}-u_{2},\Pi=\Pi_{1}-\Pi_{2}$,
where $\Pi_{1}$ and $\Pi_{2}$ are pressure functions respectively.
Then we have
\begin{equation}
 \left\{
 \begin{array}{llll}
\partial_{t}\theta+u_{1}\cdot\nabla\theta=-u\cdot\nabla\theta_{2},
\\
\theta|_{t=0}=0,
\end{array}\right.\nonumber
 \end{equation}
and
\begin{equation}
 \left\{
 \begin{array}{llll}
\partial_{t}u+u_{1}\cdot\nabla u=-u\cdot\nabla
u_{2}-\nabla\Pi+\theta e_{N},\\
{\rm div}u=0,\\
u|_{t=0}=0.
\end{array}\right.\nonumber
 \end{equation}

From Lemma 3.2, we have
$$
\|\theta(t)\|_{r-1}\leq
\int_{0}^{t}\|u\cdot\nabla\theta_{2}\|_{r-1}ds+C(r)\int_{0}^{t}\|\nabla
u_{1}\|_{L^{\infty}}\|\theta\|_{r-1}ds.
$$
Due to Lemma 2.3, we have
$$
\begin{array}{rl}
\|u\cdot\nabla\theta_{2}\|_{r-1}&\leq
\displaystyle{C(\|u\|_{L^{\infty}}\|\nabla\theta_{2}\|_{r-1}+\|\nabla\theta_{2}\|_{L^{\infty}}\|u\|_{r-1})}\nonumber\\
&\leq \displaystyle{C\|\theta_{2}\|_{r}\|u\|_{r-1}}\nonumber.
\end{array}
$$
Then we get
$$
\|\theta(t)\|_{r-1}\leq
C\int_{0}^{t}\|u\|_{r-1}\|\theta_{2}\|_{r}ds+C(r)\int_{0}^{t}\|\nabla
u_{1}\|_{L^{\infty}}\|\theta\|_{r-1}ds.
$$
Using the  Gronwall's inequality, we obtain
\begin{equation}
\|\theta(t)\|_{r-1}\leq
C\int_{0}^{T}\|\theta_{2}\|_{r}\|u\|_{r-1}ds\exp(C(r)\int_{0}^{T}\|\nabla
u_{1}\|_{L^{\infty}}ds).
\end{equation}

From Lemma 3.3, the following estimate holds true
\begin{equation}
\begin{array}{rl}
\|u\|_{r-1}&\leq\displaystyle{\|u_{0}\|_{r-1}+C(r)\int_{0}^{T}\|\nabla
u_{1}\|_{L^{\infty}}\|u\|_{r-1}+\|u\cdot\nabla
u_{2}\|_{r-1}}\nonumber\\
&+\displaystyle{\|\nabla(\Pi_{1}-\Pi_{2})\|_{r-1}+\|\theta\|_{r-1}ds},
\end{array}
\end{equation}
where
$$
\begin{array}{rl}
\nabla(\Pi_{1}-\Pi_{2})&=-\nabla\Delta^{-1}{\rm
div}(u_{1}\cdot\nabla(u_{1}-u_{2}))-\nabla\Delta^{-1}{\rm
div}((u_{1}-u_{2})\cdot\nabla u_{2})\\
&+\nabla\Delta^{-1}\partial_{N}(\theta_{1}-\theta_{2})
\end{array}
$$
i.e.
$$
\nabla\Pi=-\nabla\Delta^{-1}{\rm div}(u_{1}\cdot\nabla
u)-\nabla\Delta^{-1}{\rm div}(u\cdot\nabla
u_{2})+\nabla\Delta^{-1}\partial_{N}\theta.
$$

Using (4.18) when $1<r<2$ and using the inequality $\|\nabla
u\|_{L^{\infty}}\leq\|u\|_{r-1}$ when $r>2$, by Lemma 2.3, we get
\begin{equation}
\begin{array}{rl}
\|\nabla\Pi\|_{r-1}&\leq \displaystyle{C(\|\nabla
u_{1}\|_{L^{\infty}}\|u\|_{r-1}+ \|\nabla
u\|_{L^{\infty}}\|u_{1}\|_{r-1}}\\
&+\|\nabla u\|_{L^{\infty}}\|u_{2}\|_{r-1} +\|\nabla
u_{2}\|_{L^{\infty}}\|u\|_{r-1}+ \|\theta\|_{r-1}).
\end{array}
\end{equation}
then we have
\begin{equation}
\|\nabla\Pi\|_{r-1}\leq
\displaystyle{C(2\|u\|_{r-1}(\|u_{1}\|_{r}+\|u_{2}\|_{r})+\|\theta\|_{r-1})}.
\end{equation}
 for $1<r<2$ and $r>2$.

For $r=2$, using Lemma 2.4 and Lemma 2.3, we have
$$
\begin{array}{ll}
\|\nabla\Pi\|_{1}&\leq\displaystyle{\|u_{1}\cdot\nabla
u\|_{1}+\|u\cdot\nabla
u_{2}\|_{1}+\|\theta\|_{1}}\\
&\leq
\displaystyle{C\|u_{1}\|_{1}\|u\|_{B^{1}_{\infty,1}}+\|\theta\|_{1}}\\
&+\displaystyle{C(\|u\|_{1}\|\nabla u_{2}\|_{L^{\infty}}+\|u\|_{L^{\infty}}\|\nabla u_{2}\|_{1})}\\
&\leq\displaystyle{C(\|u_{1}\|_{1}\|u\|_{2}+\|u\|_{1}\|u_{2}\|_{2})+\|\theta\|_{1}}.
\end{array}
$$
So (4.25) is also true for $r=2$.

Putting (4.25) into (4.23), we get
$$
\|u\|_{r-1}\leq\|u_{0}\|_{r-1}+C(r)\int_{0}^{t}3\|u\|_{r-1}(\|u_{1}\|_{r}+\|u_{2}\|_{r})+2\|\theta\|_{r-1}ds.
$$
Using (4.22), we have
$$
\begin{array}{rl}
\|u\|_{r-1}&\leq\displaystyle{3C(r)\int_{0}^{T}\|u\|_{r-1}(\|u_{1}\|_{r}+\|u_{2}\|_{r})}\\
&+\displaystyle{2C(r)T\int_{0}^{T}\|\theta_{2}\|_{r}\|u\|_{r-1}ds\exp(C(r)\int_{0}^{T}\|u_{1}\|_{r}ds)}.
\end{array}
$$
In view of Gronwall's inequality, we  obtain that $u\equiv0$.
Applying (4.22) again, we directly get $\theta\equiv0$.

This completes the proof of the uniqueness in Theorem 2.1

\vspace{3mm}

\thispagestyle{empty} \setcounter{section}{5}
\setcounter{subsection}{0}
 \pagestyle{myheadings}

\setcounter{equation}{0}
 {\large \bf 5.  Blow-up criterion } \vspace{0.2cm}

\vspace{3mm}

 Now we prove Theorem 2.2, which is about the blow-up
criterion

\textbf{Proof of Theorem 2.2.} Applying  (3.4), (3.6) and the
Gronwall's inequality, we have
\begin{equation}\label{5.1}
\begin{array}{rl}
\|u\|_{r}&\leq\displaystyle{\|u_{0}\|_{r}\exp(C(r)\int_{0}^{t}\|\nabla
u\|_{L^{\infty}}ds)}\nonumber\\
&+\displaystyle{(2+2^{-r})\|\theta_{0}\|_{r}\int_{0}^{t}\exp(C(r)\int_{0}^{\tau}
\|\nabla u\|_{L^{\infty}}d\tau)ds\exp(C(r)\int_{0}^{t}\|\nabla
u\|_{L^{\infty}}ds)}.\nonumber
\end{array}
\end{equation}
Using (3.4) again and the above inequality (5.1), we get the proof
of Theorem 2.2.

\begin{center}
{\large \bf References}
\end{center}
\rm\small

 \REF{[1]}H. Abidi, T. Hmidi, On the global well-posedness for Boussinesq
system, J.Differential Equations 233(1) (2007) 199-220.

\REF{[2]}J.M.Bony, Calcul symbolique et propagation des
singularit\'{e}s pour les \'{e}quations aux d\'{e}riv\'{e}es
partielles non lin\'{e}aires, Ann.\'{E}cole Sup.14(1981)209-246.

  \REF{[3]}J. R. Cannon, E. DiBenedetto, The initial problem for the
 Boussinesq equation with data in $L^{p}$, Lecture Note in
 Mathematics,vol.771, Berlin-Heidelberg-New York, Springer,  1980, 129-144.

\REF{[4]}D. Chae, Global regularity for the 2D Boussinesq equations
with partial viscosity terms, Adv. Math. 203 (2006) 497-513.
 \REF{[5]} D. Chae, Local existence and blow up criterion for the Euler
equations in the Besov spaces, Asymptotic Analysis 38(2004)339-358.
\REF{[6]}D. Chae, S. K. Kim, H. S. Nam, Local existence and blow up
criterion of H\"{o}lder continuous solutions of the Boussinesq
equations, Nagoya Math.J.155 (1999) 55-80.
 \REF{[7]}D. Chae, H. S. Nam, Local existence and blow up criterion for the
Boussinesq equations, Proc. Roy. Soc. Edinburgh 127A(1997),935-946.
   \REF{[8]}J.Y.Chemin, Perfect Incompressible Fluids, Clarendon Press,
Oxford,1998.

\REF{[9]}J.Y.Chemin, Localization in Fourier space and Navier-Stokes
system, Lecture Notes,2005.

\REF{[10]}R.Danchin, Fourier Analysis Methods for PDF's, Lecture
Notes, 2005.

   \REF{[11]}R. Danchin, M.Paicu, Global well-posedness issues for
 the inviscid Boussinesq system with Yudovich's Type
 data, Commun. Math. Phys. 290 (2009)1-14.
 \REF{[12]} W. E, C. W. Shu, Small-scale structures in Boussinesq
 convection, Phys. Fluids 6(1) (1994), 49-58.
   \REF{[13]}J. S. Fan, Y. Zhou, A note on
regularity criterion for
 the 3D Boussinesq system with partial viscosity, Applied Mathematics
 Letters 22(2009)802-805.
 \REF{[14]}B. L. Guo,Spectral method for solving two-dimensional
 Newton-Boussinesq equation, Acta Math. Appl. Sinica, 5 (1989),
 27-50.
\REF{[15]}T. Hmidi, S. Keraani, On the global well-posedness of the
two-dimensional Boussinesq system with a zero diffusivity,
Adv.Differential Equations 12(4)(2007)461-480.

\REF{[16]} T. Y. Hou, C. Li,  Global well-posedness of the viscous
Boussinesq equations, Discrete Contin. Dyn. Syst. 12 (1)(2005)
1--12.
 \REF{[17]}N. Ishimura, H. Morimoto, Remarks on the blow up criterion 3D
 Boussinesq equations, M$^{3}$AS 9(1999)1323-1332.
     \REF{[18]}A. Majda, Introduction to PDEs and waves for the atmosphere
 and ocean, Courant Lecture Notes in
 Mathematics,AMS/CIMS, vol.9, 2003.
 \REF{[19]} J. Pedlosky, Geophysical fluid dynamics, New York,
 Springer-Verlag, 1987.
\REF{[20]} Y. Taniuchi, A note on the blow-up criterion for the
inviscid 2D Boussinesq equations,  Lecture Notes in Pure and Applied
Mathematics, 223 (2002),131-140.

 \end{document}